\documentclass[11pt]{article}

\usepackage[utf8]{inputenc}
\usepackage{graphicx}
\usepackage{labelfig}
\usepackage{amsmath,amsfonts,amssymb,amsthm,upref,bm}%
\usepackage{mathtools}
\usepackage{wasysym}
\usepackage{mathrsfs}
\usepackage{booktabs}
\usepackage{url}
\usepackage{accents}
\usepackage{cite}

\usepackage[paperwidth=210mm,paperheight=297mm,top=20mm,bottom=20mm,left=20mm,right=20mm]{geometry}

\DeclareMathAlphabet{\varmathbb}{U}{pxsyb}{m}{n}
\DeclareMathAlphabet{\mathpzc}{OT1}{pzc}{m}{it}

\renewcommand{\Re}{\operatorname{Re}}
\renewcommand{\Im}{\operatorname{Im}}

\newcommand{\MF}[1]{\mathop{\mbox{$#1$}}\nolimits}

\newcommand{\ii}{\kern0.05em\mathrm{i}\kern0.05em}
\newcommand{\E}[1]{\textrm{e}^{#1}}

\newcommand{\Zz}{\varmathbb{Z}}%
\newcommand{\Zzp}{\varmathbb{Z}_{\geqslant0}}%
\newcommand{\Zzm}{\varmathbb{Z}_{\leqslant0}}%

\newcommand{\Cb}{\varmathbb{C}}%

\newcommand{\BCzero}{\mathscr{B}_{0\infty}}%
\newcommand{\BCone}{\mathscr{B}_{1\infty}}%
\newcommand{\BCa}{\mathscr{B}_{a\infty}}

\def\Heun{H}%
\def\HeunL{\MF{\textit{Hl\/}}}%
\def\HeunS{\MF{\textit{Hs\/}}}%

\def\nst{n_{\scriptscriptstyle\star}}

\def\regsgn{\scalebox{0.5}{\bf\sun}}
\def\regsgnO{\scalebox{0.5}{\bf$\circledast$}}

\allowdisplaybreaks[3]

\begin{document}

\title{\mbox{\vrule width0mm height16mm}\bf On regularization of the Heun functions}

\author{{\large Oleg\ V.\ Motygin}}
 
\date{\it Institute for Problems in Mechanical Engineering, Russian Academy\\ of
Sciences,
 V.O., Bol'shoy pr., 61, 199178 St.\,Petersburg, Russia\\ 
\mbox{\vrule width0mm height7mm}{\rm email: \texttt{o.v.motygin@gmail.com}}}

\maketitle

\begin{abstract}

In the paper we consider the Heun functions, which are solutions of the
equation introduced by Karl Heun in 1889. The Heun functions generalize many known special functions and appear in many fields of modern physics.
Evaluation of the functions was described in \cite{Mot15}. It is based on local power series solutions near the origin, derived by the Frobenius method, and analytic continuation to the whole complex plane with branch cuts. However, exceptional cases can occur at integer values of an exponent-related parameter $\gamma$ of the equation, when one of the two local solutions should include a logarithmic term. This also means singular behavior of the Heun functions as $\gamma$ approaches the integer values.
Here we suggest a method of regularization and redefine the Heun functions in some vicinities of the integer values of $\gamma$, where the new functions depend smoothly on the parameter.

\end{abstract}

\def\texttrademark{}

\section{Introduction\mbox{\vrule width0mm height8mm}}
\label{intro0}

In the present paper we study solutions of the equation introduced in \cite{Heun} as a generalization of the
hypergeometric equation.\,\ Its general form is a Fuchsian equation with four regular singular points (see
\cite{Ronveaux1995,SlavyanovLay2000,192-Maier2007}).\,\ Now the
Heun equation (along with its confluent forms) appears in many areas of physics~--- for applications and references, see e.g.\  \cite{Hortacsu2017,BirkandanHortacsu2017,Ishkhanyan2018,Ishkhanyans2018,Mot18,Oberlack2020} and  \url{https://theheunproject.org/}.

The first software package, which was supposed to be able to evaluate the Heun functions
numerically, was Maple\texttrademark. However, the implementation is known to be not perfect.
In \cite{Mot15} an alternative procedure, based on power series expansions and analytic continuation, was suggested to define single-valued functions in the whole complex plane with branch cuts. Algorithms for numerical evaluation were worked out and a code was given in \cite{mycode}. The procedure is also applicable for computation of the multi-valued Heun functions.

Since 2020, numerical evaluation of the Heun functions is available in Wolfram Mathematica\texttrademark\ \cite{MathRelease2020}. (Various cases of confluence are also covered.) The Heun functions in Mathematica\texttrademark\ are defined in the same way as the single-valued Heun functions in \cite{Mot15}, and there is a common problem related to the so-called logarithmic values.

Evaluation of the Heun functions is based on local power series solutions near the point ${z = 0}$, derived by the Frobenius method. However, the Frobenius method generally gives two independent solutions provided that two roots of the so-called indicial equation are not separated by an integer. For the equation under consideration, the exceptional cases occur at integer values of an exponent-related parameter $\gamma$ of the equation. Then, one of the local solutions should include a logarithmic term. This also means singular behaviour of the Heun functions as functions of $\gamma$. 

In applications, it may be important to have smooth dependence of the Heun functions on the parameters.
So, in the present work we suggest a method of regularization. Comparing with the standard definitions, we redefine Heun functions in some vicinities of the integer values of $\gamma$, where the new functions depend $C^\infty$-smoothly on $\gamma$ (instead of being meromorphic).

\section{Statement and basic notations} 
\label{sect:statement}

We write Heun's equation in the standard form (see e.g.\ \cite{Erdelyi1955})
\begin{equation}
  \MF{\Heun''}(z)+\left(\frac{\gamma}{z}+\frac{\delta}{z-1}+
  \frac{\varepsilon}{z-a}\right)\MF{\Heun'}(z)+
  \frac{\alpha\beta z-q}{z(z-1)(z-a)}\MF{\Heun}(z)=0.
\label{eq:Heun}
\end{equation}
The equation has four regular singular points at $z=0$, $1$, $a$,
$\infty$ ($a \in\Cb$, $a\neq\{0,1,\infty\}$), providing the exponent-related parameters $\gamma$, $\delta$, $\varepsilon$, $\alpha$, $\beta$,  (belonging to $\Cb$) are  connected via the Fuchsian
relation
\begin{equation*}
\alpha+\beta+1=\gamma+\delta+\varepsilon.
\end{equation*}
The parameter $q\in\Cb$ is usually referred to as an accessory or auxiliary
parameter.

The Frobenius method lets one to derive local power-series solutions to
\eqref{eq:Heun} (two solutions for each of the four singular point). We study the two local solutions in a neighbourhood of the point $z = 0$. One of them corresponds to the zero root of the indicial equation, it is
called the local Heun function (see \cite{Ronveaux1995}). We denote it by
$\HeunL(a,q,\alpha,\beta,\gamma,\delta;z)$ and fix it to be equal to one for $z=0$. The second Frobenius local
solution is denoted by $\HeunS(a,q,\alpha,\beta,\gamma,\delta;z)$.

\begin{figure}[h!]
\centering\vspace{1.25mm}
 \SetLabels
 \L (0.53*0.11) {$1$}\\
 \L (0.6*0.48) {$a$}\\
 \L (0.455*0.93) {$\Im z$}\\
 \L (0.92*0.11) {$\Re z$}\\
 \L (0.19*0.28) {$\BCzero$}\\
 \L (0.71*0.28) {$\BCone$}\\
 \L (0.71*0.69) {$\BCa$}\\
 \endSetLabels
 \leavevmode\AffixLabels{\includegraphics[width=55mm]{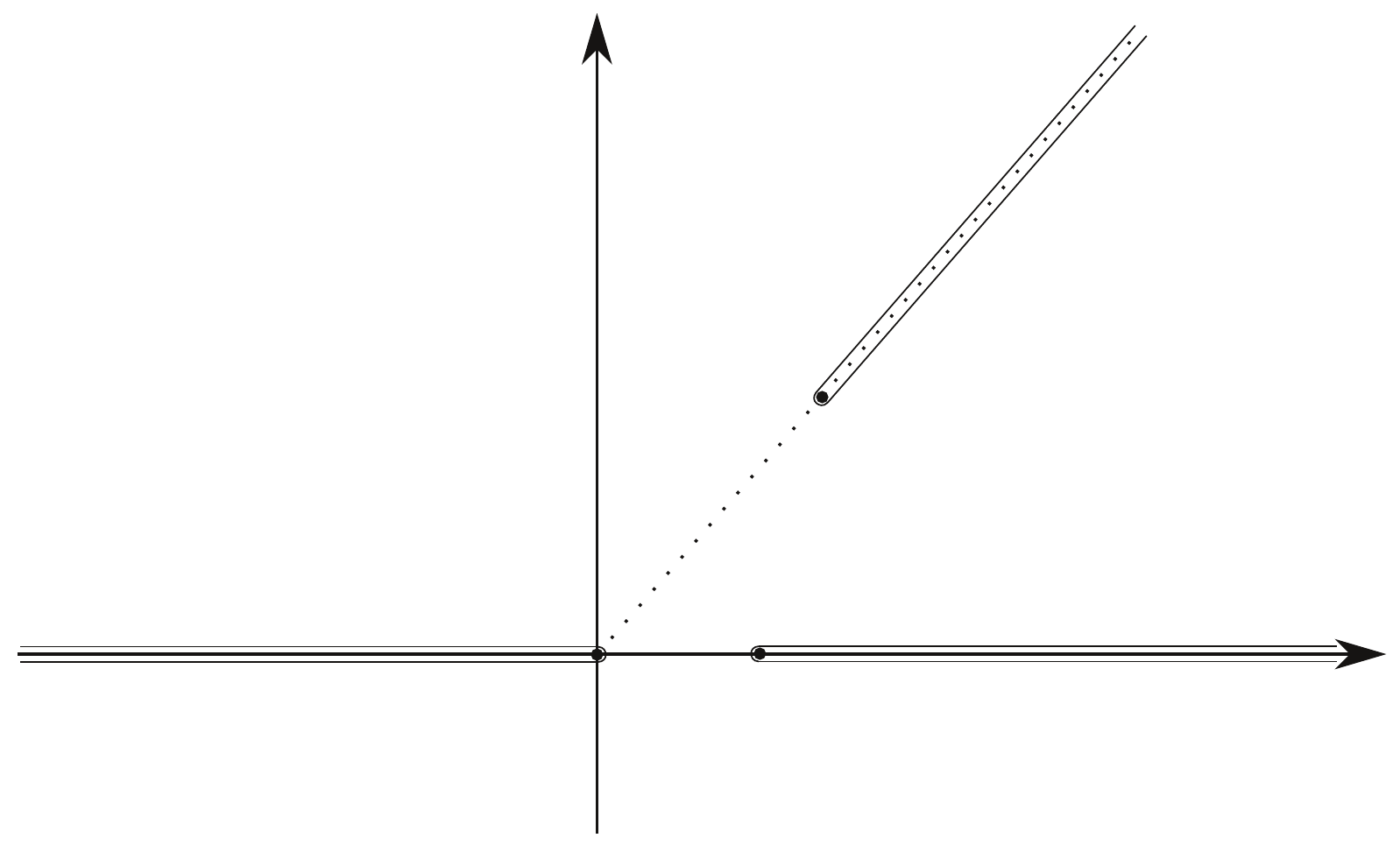}}
 \vspace{-3mm}
 \caption{Branch cuts.}
 \vspace{0mm}
 \label{bcwithoutlog}
\end{figure}
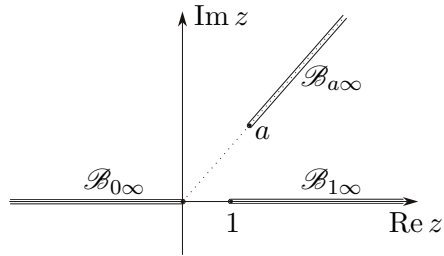

Generally, $\HeunL(a,\ldots;z)$ is a multi-valued function with
branch points at $z=1$, $a$ and $\infty$. Following \cite{Mot15}, we define a single-valued
function by fixing the branch cuts
$\BCone=(1,+\infty)$ and $\BCa=(a,\E{\ii\arg(a)}\infty)$ (shown in
Fig.~\ref{bcwithoutlog}). We also define the additional branch cut
$\BCzero=(-\infty,0)$ for the function $\HeunS(a,\ldots;z)$. All the three branch cuts are needed for the function $\HeunL(a,\ldots,z)$ for $\gamma\in\Zzm$ ($\Zzm=\{0,-1,-2,\ldots\}$ is the set of non-positive integers). 

This definition domain is star-like with respect to $z=0$ --- it is
natural as the functions $\HeunL(z)$, $\HeunS(z)$ are defined in \cite{Mot15} by analytic
continuation from a vicinity of $z=0$.

Assume first that $\gamma\not\in\Zzm$. Power series expansion of the local Heun function $\HeunL(z)$, such that $\HeunL(0)=1$, is well-known since \cite{Heun}:
\begin{equation}
\HeunL(z)=\sum_{n=0}^{\infty}b_n z^n,
\label{eq:HeunL0series}
\end{equation}
where the coefficients $b_n$ satisfy the following three-term
recurrence:
\begin{equation}
  P_n b_n=Q_n b_{n-1}+R_n b_{n-2},
\label{eq:HeunL0Andef}
\end{equation}
with the initial conditions $b_{-1}=0$, $b_0=1$. Here
\begin{equation}
\!\!\!~\begin{gathered}
P_n = a n (\gamma-1+n),\\ Q_n = q +
(n-1)\bigl[(a+1)(\gamma+n-2)+\varepsilon+a\delta\bigr],\\ R_n =
-(n-2+\alpha)(n-2+\beta).
\end{gathered}\!\!\!\!
\label{eq:PQRdef}
\end{equation}

We note that  in \eqref{eq:HeunL0Andef} $P_{\nst+1}\to0$ as $\gamma\to-\nst$, $\nst\in\Zzp:=\{0,1,2,\ldots\}$. Except for some special sets of parameters, this means singular behaviour of $\HeunL(\gamma;z)$ as a function of parameter $\gamma$. In particular, 
\[
  \frac{\partial^{\nst+1}}{\partial z^{\nst+1}}\HeunL(\gamma;z)\bigr|_{z=0}=\MF{O}\Bigl(\frac{1}{\gamma+\nst}\Bigr)\ \ \mbox{as}\ \ \gamma\to-\nst.
\]
This behaviour is illustrated below in Section~\ref{sect:oneval}.

At the same time, the local Heun function can be defined for $\gamma=-\nst$, $\nst\in\Zzp$. When $\gamma\in\Zz$, the local Frobenius solution related to the smaller exponent ($0$ or $1-\gamma$) may contain a logarithmic factor. 
Hence, for $\gamma=-\nst$, the
solution of \eqref{eq:Heun} satisfying $\HeunL(0)=1$ can be found in the form:
\begin{equation}
\HeunL(z)=\sum_{n=0}^{\infty}c_n
z^n +\bigl(C_{\nst}\log(z)+A)
\sum_{n=\nst+1}^{\infty}s_n z^n,
\label{eq:HeunLserieslog}
\end{equation}
where coefficients $c_n$ for $n=1,\ldots,\nst$
are subject to the recurrence \eqref{eq:HeunL0Andef} with  \mbox{$c_{-1}=0$}, $c_0=1$. The coefficients $s_n$ for
$n=\nst+2,\nst+3,\ldots$ are submitted to 
\eqref{eq:HeunL0Andef} with $s_{\nst}=0$, $s_{\nst+1}=1$ and
\begin{equation}
a(\nst+1)\, C_{\nst}=c_{\nst}\bigl[q-\gamma(\varepsilon+a\delta-a-1)\bigr]-
c_{\nst-1}\bigl[(1+\gamma)(2-\delta-\varepsilon)+\alpha\beta\bigr].\label{eq:Cst}
\end{equation}
We can choose $c_{\nst+1}=0$ and arbitrary coefficient $A$. The coefficients $c_n$ for $n=\nst+2,\nst+3,\ldots$
are defined by the following relationship:
\begin{equation*}
P_n c_n=Q_n c_{n-1}+R_n c_{n-2}+ C_{\nst}\bigl(S_n s_{n}+T_n s_{n-1}+U_n s_{n-2}\bigr),
\label{eq:recurrlog}
\end{equation*}
where $P_n$, $Q_n$, $R_n$ are given in \eqref{eq:PQRdef}, and
\begin{equation}
\begin{gathered}
S_n = a(1-\gamma-2n),\\\
T_n = \varepsilon+a\delta+(a+1)(\gamma+2n-3),\\
U_n = 4-2n-\alpha-\beta.
\end{gathered}
\label{eq:STUdef}
\end{equation}

So, we notice a problem: the representation \eqref{eq:HeunL0series} as a function of $\gamma$ in general case has singularities at $\gamma\in\Zzm$. At the same time, the Heun function is properly defined at these values of $\gamma$. Situation is similar for the second Heun function. We aim to redefine the Heun functions to get rid of the singularities and guarantee smooth dependence on $\gamma$. 

Defining the second function $\HeunS(z)$, we should
distinguish two situations: $\gamma=1$ and $\gamma\neq1$. For
$\gamma\neq1$, we use the following representation (see Table~2 in
\cite{192-Maier2007}, index $[0_-][1_+][a_+][\infty_-]$):
\begin{equation}
 \!\HeunS(a,q,\alpha,\beta,\gamma,\delta;z) \\{}= z^{1-\gamma}\HeunL(a,q',\alpha',\beta',\gamma',\delta;z),
\label{eq:HeunSgammaneq1}
\end{equation}
where $q'=q-(\gamma-1)(\varepsilon+a\delta)$, $\alpha'=\alpha-\gamma+1$,
$\beta'=\beta-\gamma+1$, $\gamma'=2-\gamma$.

The representation \eqref{eq:HeunSgammaneq1} shows that the second Heun function has singularities at $\gamma=2,3,4,\ldots$
On the other hand, for $\gamma\neq1$ we do not have to regularize the function $\HeunS$ separately, providing we have a regularized version of the function $\HeunL$.

It is important to note that for $\gamma\in\Zzm$, by \eqref{eq:HeunL0series} and \eqref{eq:HeunSgammaneq1}, we have
\begin{equation}
  \HeunS(z)=\sum_{n=\nst+1}^{\infty}s_n z^n,
\label{eq:HeunSsergammanonpositive}
\end{equation}
and can rewrite \eqref{eq:HeunLserieslog} as follows:
\begin{equation}
\HeunL(z)=\sum_{n=0}^{\infty}c_n
z^n +\bigl(C_{\nst}\log(z)+A)
\HeunS(z),
\label{eq:HeunLserieslog'}
\end{equation}

It is easy to note that for $\gamma=1$ the representation \eqref{eq:HeunSgammaneq1} of $\HeunS(z)$  coincides with $\HeunL(z)$. This is another type of degeneration of the pair of solutions to the Heun equation. However, as shown in \cite{Mot15}, at $\gamma=1$ we also have a proper definition of the function $\HeunS$, including a logarithmic term:
\begin{equation}
\HeunS(z)=\sum_{n=1}^{\infty}d_n z^n+(\log(z)+B)\HeunL(z),
\label{eq:HeunSgamma=1ser}
\end{equation}
where $B$ is an arbitrary coefficient and $d_n$ satisfy
\begin{equation*}
P_n d_n=Q_n d_{n-1}+R_n d_{n-2}+ S_n t_{n}+T_n t_{n-1}+U_n t_{n-2},
\end{equation*}
with $d_{-1}=d_0=0$ and coefficients are defined by \eqref{eq:PQRdef} and \eqref{eq:STUdef}.

Further we will redefine the Heun functions $\HeunL$ and $\HeunS$ to avoid both types of degeneration at the integer values of $\gamma$.

\section{Regularization of the function $\bm{\HeunL(\gamma;z)}$ at non-positive integer $\bm\gamma$}

In \cite{Mot18} (where the confluent Heun functions were studied), it was suggested to define a regularized function dividing by a function having poles at \mbox{$\gamma\in\Zzm$}. It is natural to write
\[
\HeunL_{\regsgnO}(\gamma;z):=
\frac{1}{\Gamma(\gamma)}\HeunL(\gamma;z),
\]
where $\Gamma(\cdot)$ is the gamma function.

However, this definition has obvious shortcomings (as was partly noted in \cite{Mot18}). Consider the limit of $\HeunL_{\regsgnO}(\gamma;z)$ as $\gamma\to-\nst$, $\nst\in\Zzp$, which, on the other hand, can be obtained from the general form of solution \eqref{eq:HeunLserieslog'} by a choice of coefficients. In view of \eqref{eq:HeunL0series}, it is easy to note that $c_n=0$ for $n=0,1,\ldots,\nst$. Thus, by \eqref{eq:Cst} $C_{\nst}=0$ and $\HeunL_{\regsgnO}(\gamma;z)$ coincides with $\mathrm{const}\cdot\HeunS(\gamma;z)$ when $\gamma=-\nst$. Besides, it possible that, for some particular sets of parameters, the right-hand side of \eqref{eq:HeunL0Andef}  for $n=1+\nst$ tends to zero as $\gamma\to-\nst$. This would mean that at this value of $\gamma$ $\HeunL(\gamma;z)$ does not have the singularity and, hence, $\HeunL_{\regsgnO}(\gamma;z)\equiv0$. Along with the degeneration, behaviour of the regularized function is strongly affected by the gamma function for $|\gamma|\gg1$. 

So, in the present note we propose a new approach to regularization based on combining  $\HeunL(\gamma;z)$ and $\HeunS(\gamma;z)$. First, consider values of $\gamma$ close to a non-positive integer $\gamma=-\nst$, $\nst\in\Zzp$. We suggest the following algorithm:\\[1.5mm]
$\bullet$ fix $\gamma=-\nst$ and compute coefficients $c_n$ for $n=1,\ldots,\nst$,
using the recurrence relation \eqref{eq:HeunL0Andef} with the
initial conditions $c_{-1}=0$, $c_0=1$ (cf.\ \eqref{eq:HeunLserieslog});\\[1.5mm]
$\bullet$  compute the value $$K_{\nst}=\frac{Q_{\nst+1} c_{\nst}+R_{\nst+1} c_{\nst-1}}{a (\nst+1)};$$\\
$\bullet$ define a regularized Heun function in a vicinity of $\gamma=-\nst$ as follows:
\begin{equation}
\HeunL(\gamma;z)-\frac{K_{\nst}}{\gamma+\nst}\HeunS(\gamma;z).
\label{eq:Hreg0}
\end{equation}

Considering the limit $\gamma\to-\nst$, we note that in the power series of \eqref{eq:Hreg0}, the coefficients of $z^n$, \mbox{$n=0,1,\ldots,\nst$},  do not cancel --- see \eqref{eq:HeunSgammaneq1} and \eqref{eq:HeunSsergammanonpositive}. So, unlike $\HeunL_{\regsgnO}(\gamma;z)$, as $\gamma\to-\nst$ we come to a function of $z$ linearly independent to $\HeunS(z)$, i.e.\ to the solution \eqref{eq:HeunLserieslog} with some constant $A$.

\begin{figure}[b!]
\centering\vspace{1.25mm}
 \SetLabels
 \L (-0.14*0.87) $\rho(r)$\\
 \L (0.87*-0.055) $r$\\
 \endSetLabels
 \leavevmode\AffixLabels{\includegraphics[width=39mm]{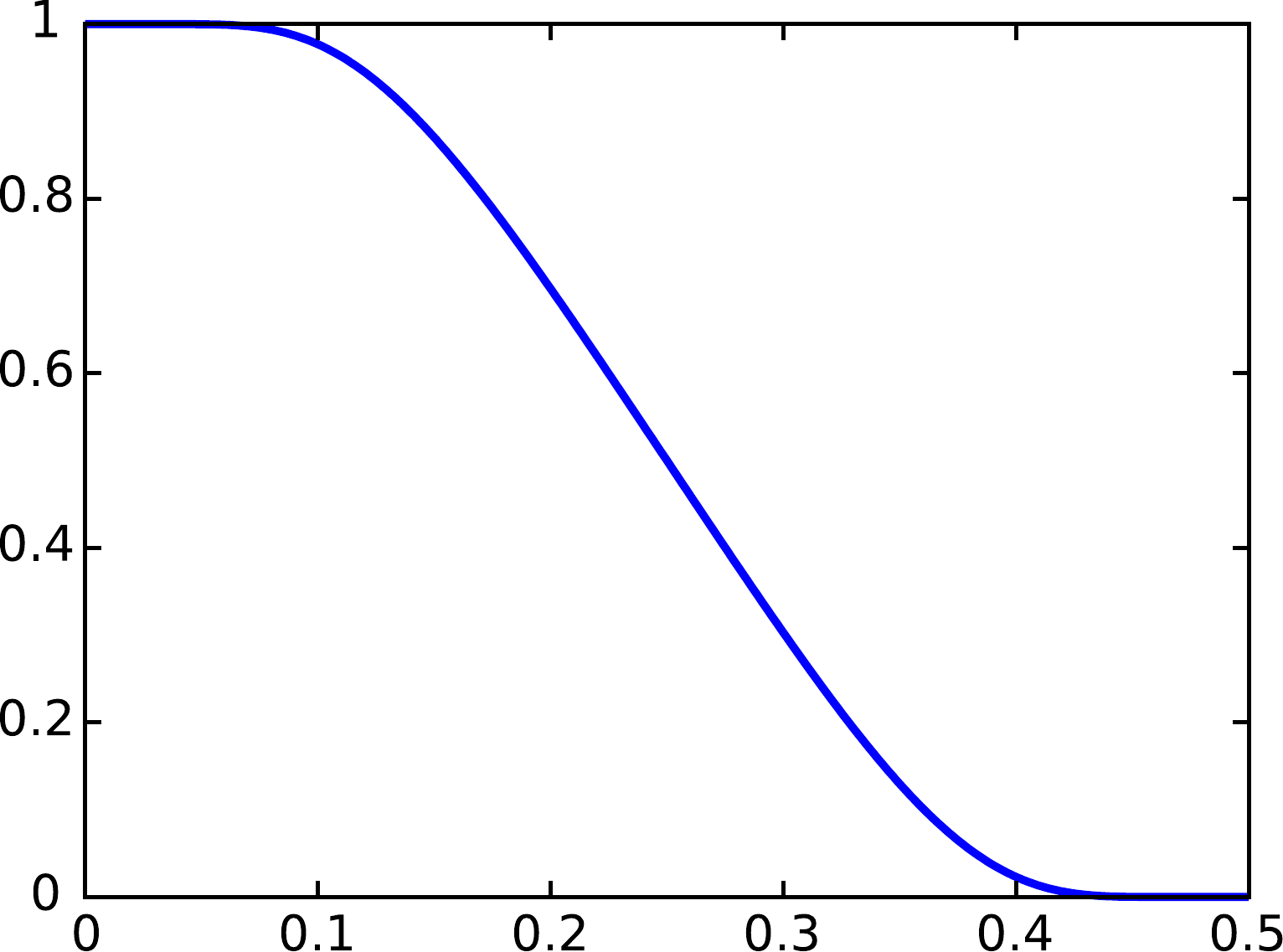}}
 \caption{The cut-off function $\rho(r)$.}
\label{fig:3}
\end{figure}

Now we would like to define a function $\HeunL_{\regsgn}(\gamma;z)$ that is regularized at all non-positive integer values of $\gamma$. We will do it by local correction of $\HeunL(\gamma;z)$ in vicinities of $\gamma\in\Zzm$. For $\gamma\in\Cb$, we define
\begin{equation}
\HeunL_{\regsgn}(\gamma;z)=\HeunL(\gamma;z),\ \
\mbox{when}\ \ \operatorname{dist}\bigl(\gamma,\Zzm\bigr)\geq1/2,
\label{eq:Hlreg1}
\end{equation}
and\begin{samepage}
\begin{equation}
\HeunL_{\regsgn}(\gamma;z)=\HeunL(\gamma;z)-
\frac{K_{\nst}}{\gamma+\nst}\rho(|\gamma+\nst|)\HeunS(\gamma;z),
\label{eq:Hlreg2}
\end{equation}
when $|\gamma+\nst|<1/2$ for some $\nst\in\Zzp$.\end{samepage}
Here $\rho(r)\in C^\infty(\mathbb{R})$ is a cutoff function. Further, we choose the following function shown in Fig.~\ref{fig:3}:
\begin{equation}
\rho(r)=\begin{cases}
1, & r\leq0,\\
\frac{\exp\{1/(2r)+1/(2r-1)\}}{1+\exp\{1/(2r)+1/(2r-1)\}}, & 0<r<1/2,\\
0, & r\geq1/2.
\end{cases}
\label{eq:rhodef}
\end{equation}

\begin{figure}[t!]
\centering
\leavevmode\AffixLabels{\includegraphics[width=62mm]{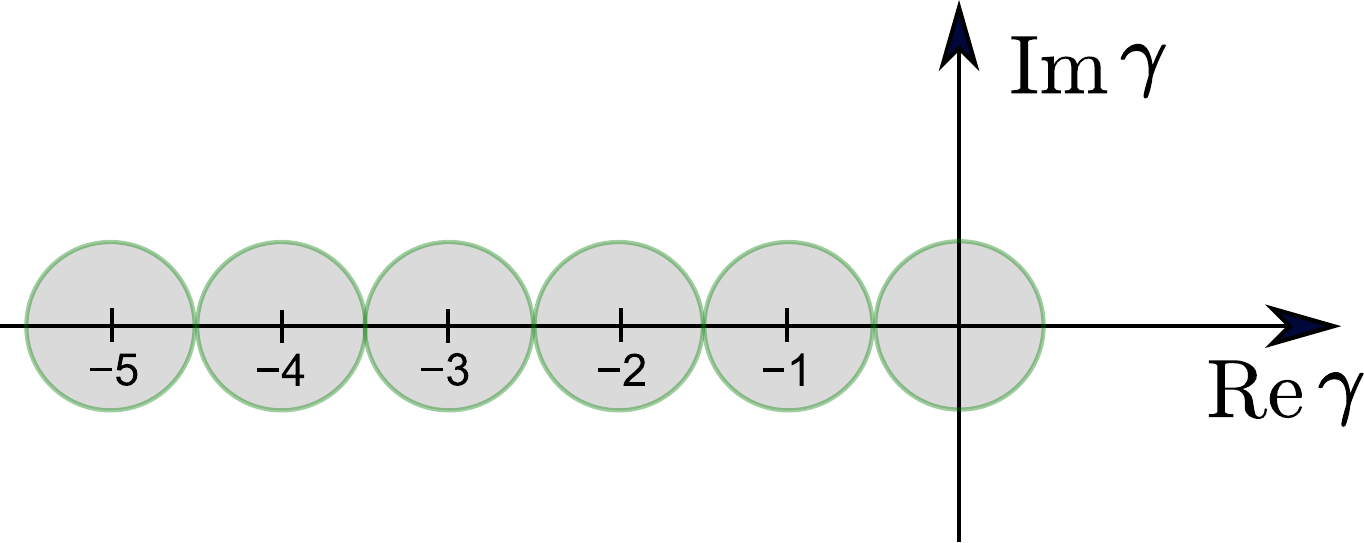}}
\caption{Vicinities of $\gamma=0,-1,-2,-3,\ldots$, where the representation \eqref{eq:Hlreg2} is used and the additional branch cut $\BCzero$ for $\HeunL_{\protect\regsgn}(\gamma;z)$ is needed.}
\label{fig:4}
\end{figure}

We remind that, generally, definition of the single-valued $\HeunL(z)$ involves two branch cuts $\BCone$ and $\BCa$ when $\gamma\in\Cb\setminus\Zzm$. For $\gamma\in\Zzm$ (logarithmic case), generally, definition of $\HeunL(z)$ demands also the branch cut $\BCzero$ (see Fig.~\ref{bcwithoutlog}). Due to the involvement of $\HeunS(z)$ into \eqref{eq:Hlreg2}, definition of the regularized Heun function $\HeunL_{\protect\regsgn}(z)$ needs the additional branch cut $\BCzero$ for the set of $\gamma$ shown in Fig.~\ref{fig:4} $\{\gamma:\operatorname{dist}\bigl(\gamma,\Zzm\bigr)<1/2\}$. We also note that $\HeunL_{\regsgn}(\gamma;z)$ is analytic in \mbox{$\gamma\in\Cb$} except these vicinities of $0,-1,-2,\ldots$, where the function is $C^\infty$-smooth.

\section{Regularization of the function $\bm{\HeunS(\gamma;z)}$ at positive integer $\bm\gamma$}

Using \eqref{eq:HeunSgammaneq1} with the function defined by \eqref{eq:Hlreg1}, \eqref{eq:Hlreg2}, we introduce
\begin{equation*}
\accentset{\circ}{\HeunS}_{\regsgn}(a,q,\alpha,\beta,\gamma,\delta;z) = z^{1-\gamma}\HeunL_{\regsgn}\bigl(a,q-(\gamma-1)(\varepsilon+a\delta),
 \alpha-\gamma+1,\beta-\gamma+1,2-\gamma,\delta;z\bigr).
\end{equation*}
This function is analytic in $\gamma$ in $\Cb$ except the vicinities of $2,3,4,\ldots$, where the function is $C^\infty$-smooth.

However, the system of Heun functions $\{\HeunL_{\regsgn}(\gamma;z), \accentset{\circ}{\HeunS}_{\regsgn}(\gamma;z)\}$ degenerates at $\gamma=1$:
$$\accentset{\circ}{\HeunS}_{\regsgn}(z) =\HeunL_{\regsgn}(z).$$
In a vicinity of $\gamma=1$, we suggest to specify the second Heun function as follows:
\begin{equation}
\frac{1}{1-\gamma}\bigl(\HeunS(\gamma;z)-\HeunL(\gamma;z)\bigr).
\label{eq:reg1}
\end{equation}
Taking into account \eqref{eq:HeunSgammaneq1} and the asymptotic representation $z^{1-\gamma}=1+(1-\gamma)\log(z)+\MF{O}\bigl((\gamma - 1)^2\bigr)$ as $\gamma\to1$, we find the limit of \eqref{eq:reg1}:
\begin{equation*}
\log(z)\HeunL(a,q,\alpha,\beta,\gamma,\delta;z) -
\frac{\textrm{d}~}{\textrm{d}\gamma}\HeunL\bigl(a,q-(\gamma-1)(\varepsilon+a\delta),
 \alpha-\gamma+1,\beta-\gamma+1,2-\gamma,\delta;z\bigr)\bigr|_{\gamma=1}.
\end{equation*}
In view of analytic dependence of $\HeunL(\gamma;z)$ on $\gamma$ in a vicinity of $\gamma=1$, it is easy to note that the latter formula corresponds to \eqref{eq:HeunSgamma=1ser} with some $B$.

Finally, we define
\begin{equation*}
\HeunS_{\regsgn}(\gamma;z)=
\accentset{\circ}{\HeunS}_{\regsgn}(\gamma;z),\ \ \mbox{when}\ \ |\gamma-1|\geq1/2,
\label{eq:HeunSgamma1reg1}
\end{equation*}
 and
\begin{equation}
\HeunS_{\regsgn}(\gamma;z)=
\rho(|\gamma-1|)\frac{1}{1-\gamma}\bigl(\accentset{\circ}{\HeunS}_{\regsgn}(\gamma;z)-\HeunL(\gamma;z)\bigr)
+\bigl(1-\rho(|\gamma-1|)\bigr)\accentset{\circ}{\HeunS}_{\regsgn}(\gamma;z),
\label{eq:HeunSgamma1reg2}
\end{equation}
when $ |\gamma-1|<1/2$.

\begin{figure*}[t!]\centering\vspace{1.25mm}
 \SetLabels
 \L (0.12*0.945) $\Re\HeunL_{\regsgn}(\gamma; z)$\\
 \L (0.04*0.24) $\Im \gamma$\\
 \L (0.7*0.13) $\Re \gamma$\\
 \endSetLabels
 \leavevmode\AffixLabels{\includegraphics[width=72mm]{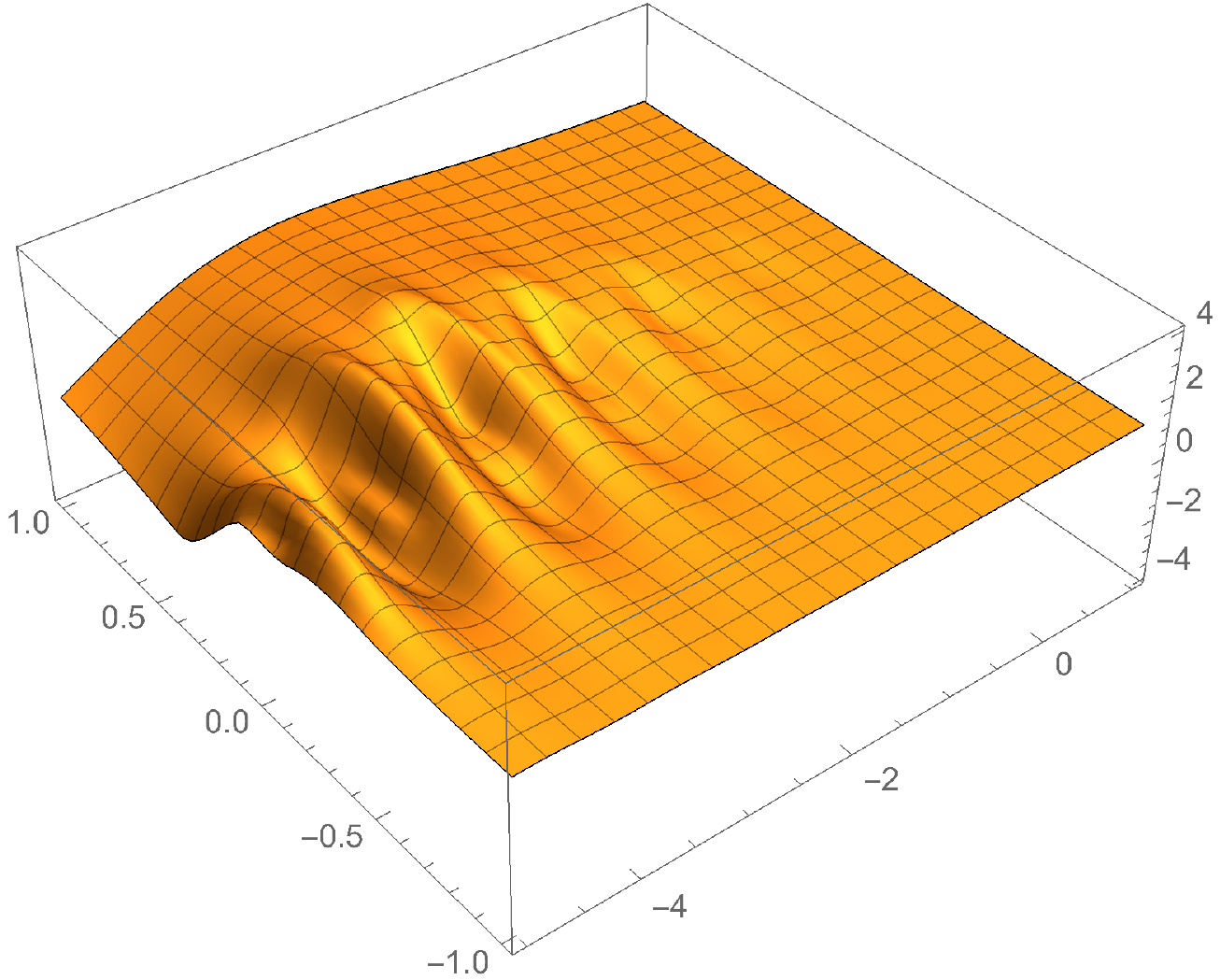}}
\kern14mm
 \SetLabels
 \L (0.14*0.945) $\Re\HeunL(\gamma; z)$\\
 \L (0.04*0.24) $\Im \gamma$\\
 \L (0.7*0.13) $\Re \gamma$\\
 \endSetLabels
 \leavevmode\AffixLabels{\includegraphics[width=72mm]{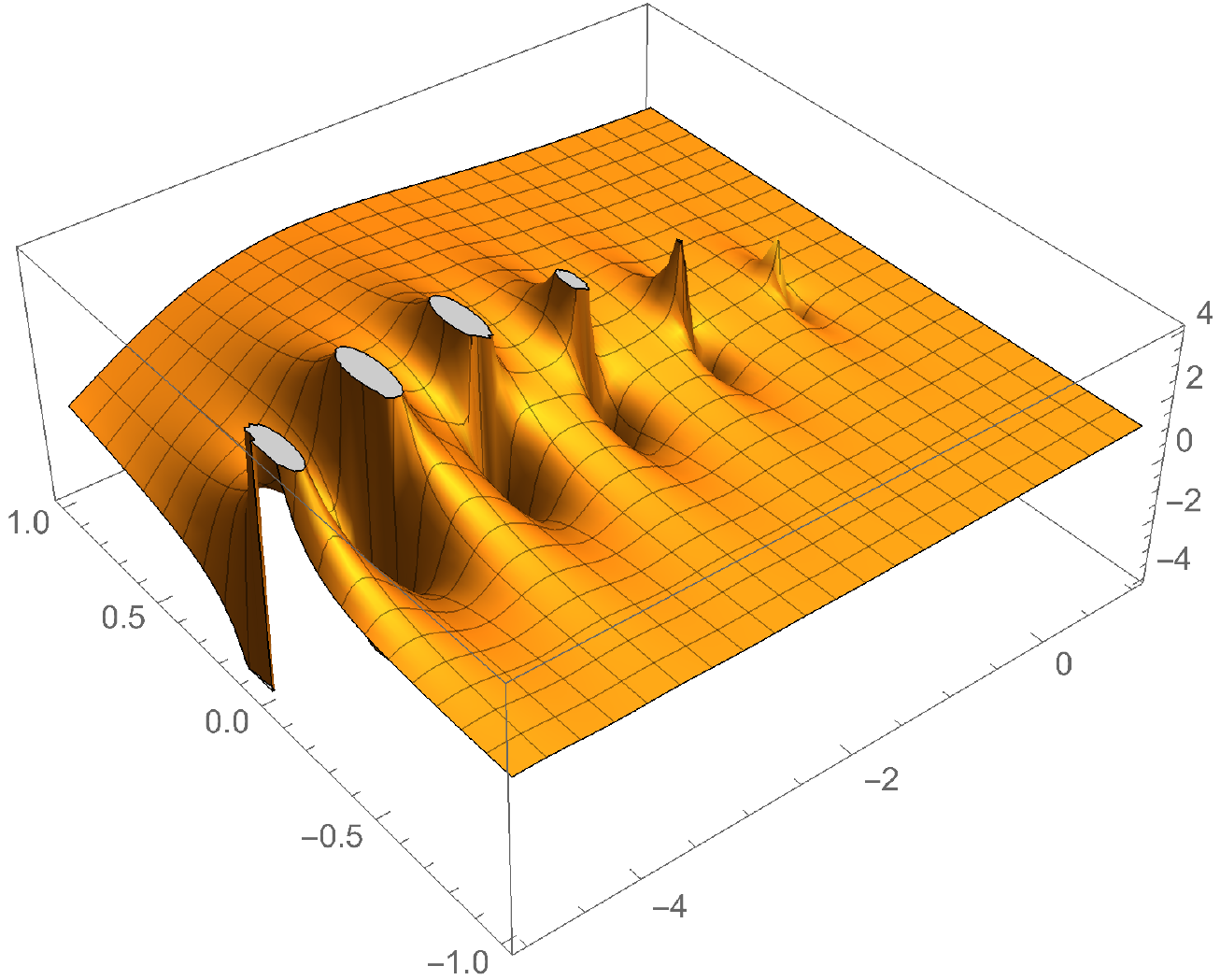}}
\caption{Dependence of $\Re\HeunL_{\protect\regsgn}(a,q,\alpha,\beta,\gamma,\delta;z)$ and $\Re\HeunL(a,q,\alpha,\beta,\gamma,\delta;z)$ on complex $\gamma$ for fixed $a=1 + \ii$, $q=0.3$, $\alpha=1.4 + 0.9 \ii$, $\beta=1.1$, $\delta=6.7$, $z=\ii$.}
\label{fig:num2}
\end{figure*}

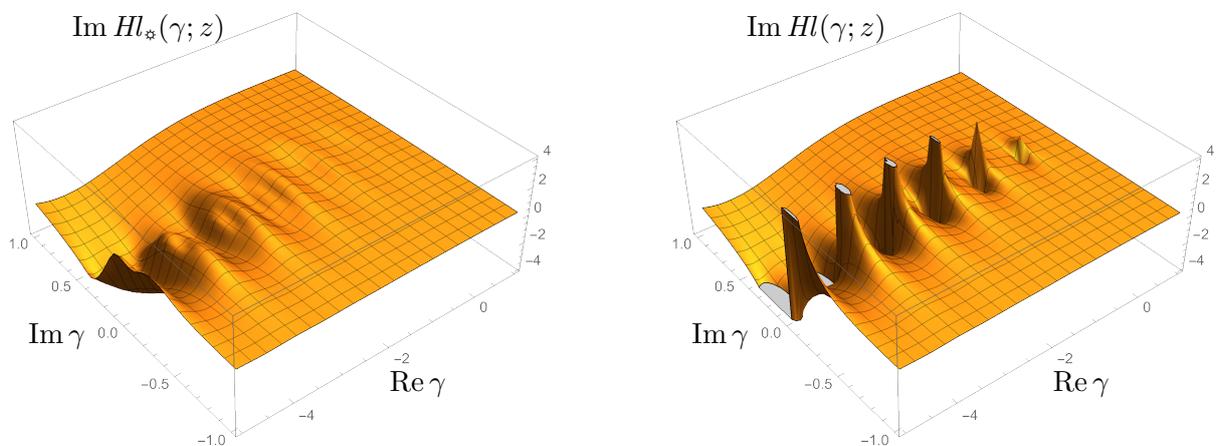
\begin{figure*}[t!]\centering\vspace{1.25mm}
 \SetLabels
 \L (0.12*0.945) $\Im\HeunL_{\regsgn}(\gamma; z)$\\
 \L (0.04*0.24) $\Im \gamma$\\
 \L (0.7*0.13) $\Re \gamma$\\
 \endSetLabels
 \leavevmode\AffixLabels{\includegraphics[width=72mm]{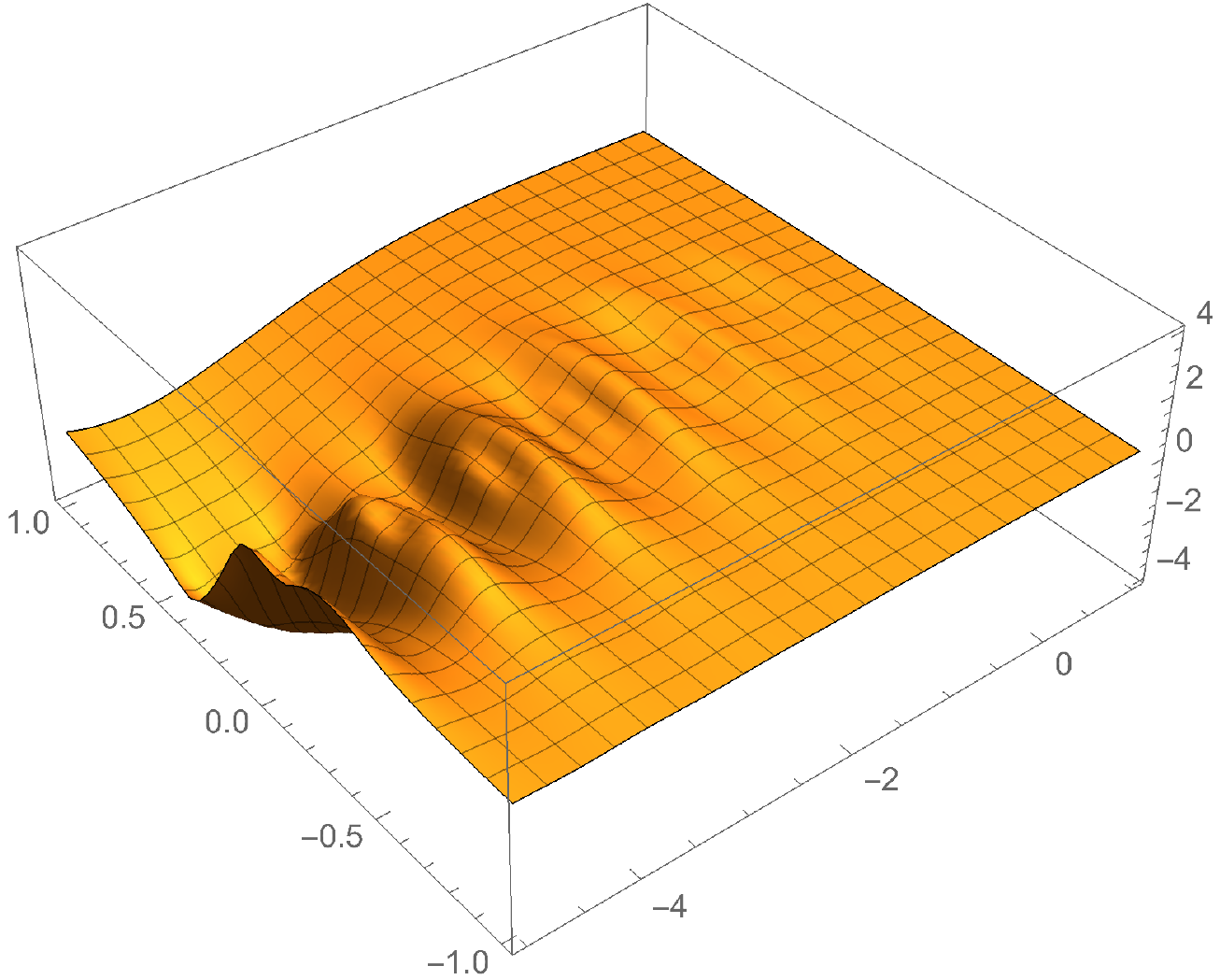}}
\kern14mm
 \SetLabels
 \L (0.14*0.945) $\Im\HeunL(\gamma; z)$\\
 \L (0.04*0.24) $\Im \gamma$\\
 \L (0.7*0.13) $\Re \gamma$\\
 \endSetLabels
 \leavevmode\AffixLabels{\includegraphics[width=72mm]{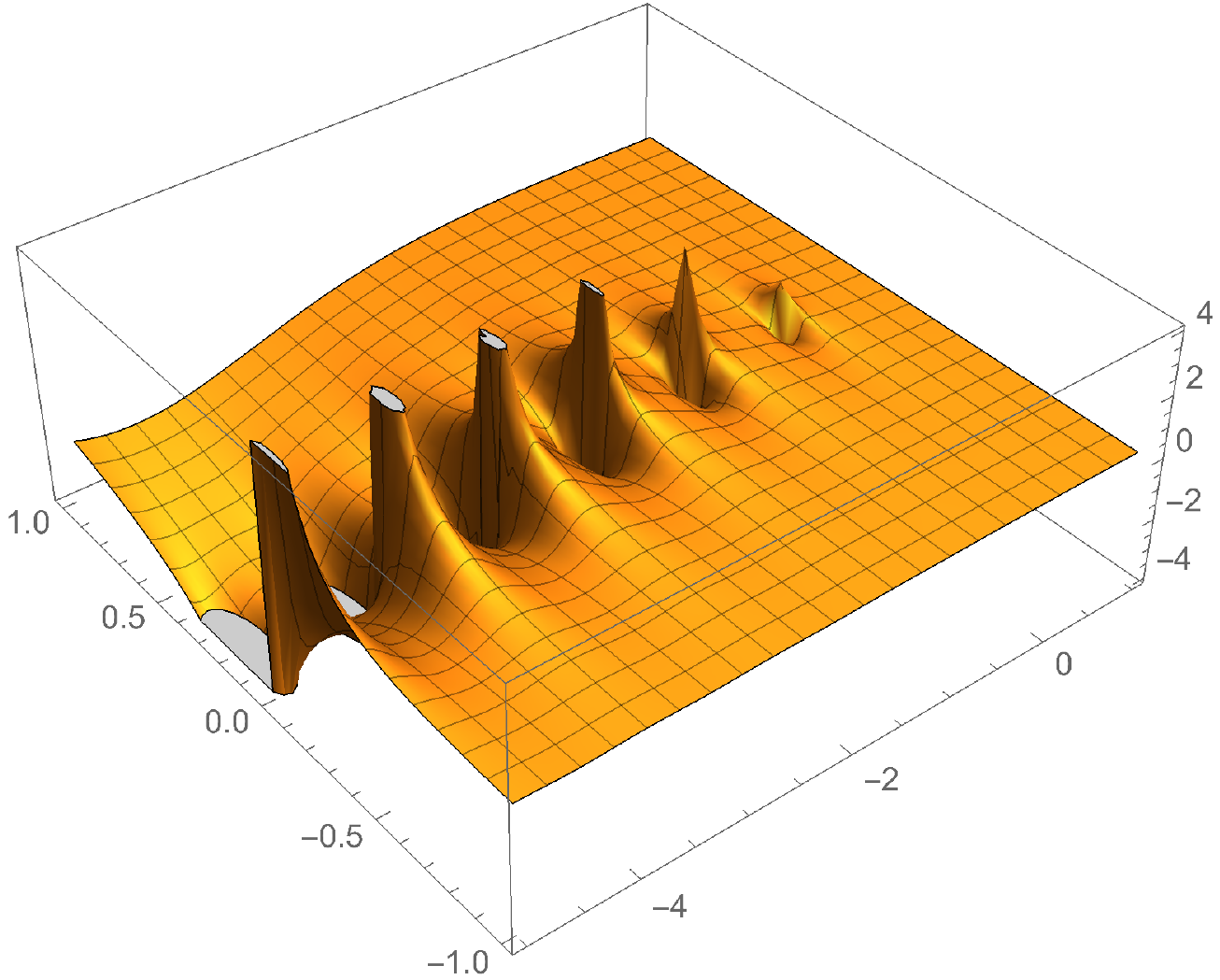}}
\caption{Dependence of $\Im\HeunL_{\protect\regsgn}(a,q,\alpha,\beta,\gamma,\delta;z)$ and $\Im\HeunL(a,q,\alpha,\beta,\gamma,\delta;z)$ on complex $\gamma$ for fixed $a=1 + \ii$, $q=0.3$, $\alpha=1.4 + 0.9 \ii$, $\beta=1.1$, $\delta=6.7$, $z=\ii$.}
\label{fig:num3}
\end{figure*}

\section{Evaluation of the regularized Heun functions}
\label{sect:oneval}

We have the asymptotics
\begin{equation}
\rho(r)=1+\MF{O}(\exp(-1/2r))\ \ \mbox{as}\ \ r\to0^+,
\label{eq:rhoest}
\end{equation}
for the function $\rho(r)$ given by \eqref{eq:rhodef}.
So, in view of \eqref{eq:Hlreg2}, for $\gamma$ close to an integer $m\leq0$, $\HeunL_{\regsgn}(\gamma;z)$ is close to the analytic function \eqref{eq:Hreg0}. For its evaluation, it is convenient to use Cauchy formula. Hence, we can exploit the following:\begin{samepage}
\begin{equation*}
\HeunL_{\regsgn}(\gamma;z)\rightarrow\frac{1}{2\pi\mathrm{i}}
\oint_{\Upsilon}\frac{\HeunL_{\regsgn}(\gamma';z)}{\gamma'-\gamma}\textrm{d}\gamma',\ \
\mbox{as}\ \ \max_{\varsigma\in\Upsilon}\{|\varsigma-m|\}\to0,
\end{equation*}
where $\Upsilon$ is a small contour enclosing $\gamma$ and $m$ and passed counter-clockwise.\end{samepage} In numerical evaluation of the function $\HeunL_{\regsgn}(\gamma;z)$, it convenient that size of $\Upsilon$ can be reasonably large (in view of \eqref{eq:rhoest}) comparing to the needed accuracy of computations. Analogous method can be applied to compute $\HeunS_{\regsgn}(\gamma;z)$ near $\gamma=1$ as defined by the formula \eqref{eq:HeunSgamma1reg2}.

Computation of the Heun functions $\HeunL_{\regsgn}(\gamma;z)$, $\HeunS_{\regsgn}(\gamma;z)$ in a vicinity of $z=0$ can be performed using the series \eqref{eq:HeunL0series}. A procedure of analytic continuation to the Mittag-Leffler star (shown in Fig.~\ref{bcwithoutlog}), suggested in \cite{Mot15}, is based on  the power
series expansion of a solution $\Heun(z)$ to the equation
\eqref{eq:Heun} having prescribed values 
\mbox{$\Heun(z_0)$},
\mbox{$\Heun'(z_0)$}
at an arbitrary finite point $z=z_0$, $z_0\not\in\{0,1,a\}$. 
First we find the value of the function at $z_*=\varrho\E{\ii\arg(z)}$ for sufficiently small $\varrho>0$ and then use the expansion consequently in a system of overlapping circular elements along a path ending at the  value $z$. The procedure of \cite{Mot15} is applicable to the regularized function without any restriction.

In Figs.~\ref{fig:num2}, \ref{fig:num3}, we present results of numerical evaluation of the functions 
$\HeunL_{\regsgn}(\gamma;z)$, $\HeunS_{\regsgn}(\gamma;z)$ and compare these regularized function with the original ones.

\end{document}